\documentclass[12pt]{article}
\usepackage{latexsym, amssymb, amscd, amsfonts, epsfig, graphicx, colordvi,verbatim,ifpdf}
\usepackage{amsmath}
\usepackage{epsfig,cite, psfrag,eepic,color}
\usepackage{amscd,graphics}
\usepackage{graphicx}
\usepackage{fancybox}
\usepackage{float}
\usepackage{pstricks}
\usepackage{mathrsfs}
\usepackage{cases}

\newtheorem{thm}{Theorem}[section]
\newtheorem{pro}[thm]{Proposition}

\newtheorem{defi}[thm]{Definition}
\newtheorem{lem}[thm]{Lemma}

\def\qed{\nopagebreak\hfill{\rule{4pt}{7pt}}
\medbreak}

\numberwithin{equation}{section}

\def\qed{\nopagebreak\hfill{\rule{4pt}{7pt}}
\medbreak}

\newlength{\boxedparwidth}
  {\begin{center} \begin{tabular}{|@{\hspace{.315in}}c@{\hspace{.15in}}|}
                  \hline \\ \begin{minipage}[t]{\boxedparwidth}
                  \setlength{\parindent}{.25in}}%
  {\end{minipage} \\ \\ \hline \end{tabular} \end{center}}

\begin{document}
\begin{center}
{\large \bf Ramanujan-type Congruences for  $\ell$-Regular Partitions \\[7pt] Modulo $3, 5, 11$ and $13$ }
\end{center}

\begin{center}
{Hai-Tao Jin}$^{1}$ and {Li Zhang}$^{2}$

\vskip 2mm

 $^{1}$School of Science,
   Tianjin University of Technology and Education \\[5pt]
   Tianjin 300222, P. R. China \\[6pt]
 $^{2}$Center for Combinatorics, LPMC-TJKLC\\
   Nankai University, Tianjin 300071, P. R. China\\

 \vskip 2mm

 Emails: $^{a}$jinht1006@tute.edu.cn,   $^{b}$zhangli427@mail.nankai.edu.cn
\end{center}

\vskip 6mm \noindent {\bf Abstract.} Let $b_\ell(n)$ be the number of $\ell$-regular partitions of $n$. Recently, Hou et al established several infinite families of congruences for $b_\ell(n)$ modulo $m$, where $(\ell,m)=(3,3),(6,3),(5,5),(10,5)$ and $(7,7)$. In this paper, by the vanishing property given by Hou et al, we show an infinite family of congruence for
$b_{11}(n)$ modulo $11$. Moreover, for $\ell= 3, 13$ and $25$, we obtain three infinite families of congruences for $b_{\ell}(n)$ modulo $3, 5$ and $13$  by the theory of Hecke eigenforms.



\section{Introduction}
A partition of a nonnegative integer $n$ is a nonincreasing sequence of positive integers whose sum is $n$.
For a positive integer $\ell$, we say a partition of $n$ is an $\ell$-regular partition provided none of its parts is divisible by $\ell$. Denote the number of $\ell$-regular partitions of $n$ by $b_\ell(n)$. For convenience, we set $b_\ell(0)=1$ and $b_\ell(n)=0$ if $n<0$. The generating function of $b_\ell(n)$ is given by
\[
B_\ell(q) = \sum_{n=0}^\infty b_\ell(n)q^n = \frac{(q^\ell; q^\ell)_\infty}{(q; q)_\infty},
\]
where
\[
(q;q)_\infty = \prod_{i=1}^\infty (1-q^i)
\]
is the standard notation in $q$-series.

Recently,
%distribution and the divisibility of $b_\ell(n)$ have been widely studied, see, for example \cite{Alladi-1997,Lovejoy-2001,Gordon-Ono-1997,Ono-Penniston-2000}.
the arithmetic properties of $b_\ell(n)$ have been widely studied. For example, Andrews, Hirschhorn and Sellers \cite{Andrews-Hirschhorn-Sellers-2010} derived some infinite families of congruences for $b_4(n)$ modulo $m$ with $m \in \{2,3,4,6,12\}$. Dandurand and Penniston\cite{Dandurand-Pennistion-2009} gave precise conditions of $n$ such that
$b_\ell(n)\equiv 0 \pmod{\ell}$ with $\ell\in\{5, 7, 11\}$. Cui and Gu \cite{Cui-Gu-2013} derived  congruences for $b_\ell(n)$ modulo $2$ with $\ell\in \{2,4,5,8,13,16\}$ by employing the $p$-dissection formulas of Ramanujan's theta functions $\psi(q)$ and $f(-q)$. Keith \cite{Keith-2014} studied the congruences for $b_9(n)$ modulo $3$. Furcy and Penniston \cite{Furcy-Penniston-2012} obtained congruences for $b_\ell(n)$ modulo $3$ with $\ell \in \{4,7,13,19,25,34,37,43,49\}$ by using the theory of modular forms. Hou et al \cite{Hou-Sun-Zhang-2015} proved several infinite families of congruences  for $b_\ell(n)$ modulo $3$, $5$ and $7$ by applying binary quadratic form approach.

The objective of this paper is to derive several infinite families of congruences for $\ell$-regular partition modulo $3$, $5$, $11$ and $13$. Our approaches are based on the theory of Hecke eigenforms and  binary quadratic forms due to \cite{Hou-Sun-Zhang-2015}. The main results of this paper are stated as follows.

\begin{thm}\label{b3-p-thm}
Let $\alpha, n$ be nonnegative integers and  $p_i\geq 5$ be primes such that $p_i\not\equiv 1 \pmod{12}$. Then for any integer $j \not\equiv 0 \pmod{p_{\alpha+1}}$,
we have
\begin{align}\label{b3-p-cong}
 b_3\left(p_1^2\cdots p_{\alpha+1}^2n+\frac{p_1^2\cdots p_{\alpha}^2p_{\alpha +1}(12j+p_{\alpha +1})-1}{12}\right)\equiv 0 \pmod{3}.
\end{align}
\end{thm}

\begin{thm}\label{b11-p-thm}
Let $\alpha, n$ be nonnegative integers and  $p_i\geq 5$ be primes such that $p_i\equiv 3 \pmod{4}$. Then for any integer $j \not\equiv 0 \pmod{p_{\alpha+1}}$,
we have
\begin{align}\label{b11-p-cong}
 b_{11}\left(p_1^2\cdots p_{\alpha+1}^2n+\frac{p_1^2\cdots p_{\alpha}^2p_{\alpha +1}(12j+5p_{\alpha +1})-5}{12}\right)\equiv 0 \pmod{11}.
\end{align}
\end{thm}

\begin{thm}\label{b13-p-thm}
Let $\alpha, n$ be nonnegative integers and  $p_i$ be odd primes such that $b_{13}(\frac{p_i-1}{2})\equiv 0 \pmod{13}$. Then for any integer $j \not\equiv 0 \pmod{p_{\alpha+1}}$ , we have
\begin{align}\label{b13-p-cong}
 b_{13}\left(p_1^2\cdots p_{\alpha+1}^2n+\frac{p_1^2\cdots p_{\alpha}^2p_{\alpha +1}(2j+p_{\alpha +1})-1}{2}\right)\equiv 0 \pmod{13}.
\end{align}
\end{thm}

\begin{thm}\label{b25-p-thm}
Let $\alpha, n$ be nonnegative integers and  $p_i$ be primes such that $b_{25}(p_i-1)\equiv 0 \pmod{5}$. The for any integer $j \not\equiv 0 \pmod{p_{\alpha+1}}$,
we have
\begin{align}\label{b25-p-cong}
 b_{25}\left(p_1^2\cdots p_{\alpha+1}^2n+p_1^2\cdots p_{\alpha}^2p_{\alpha +1}(j+p_{\alpha +1})-1\right)\equiv 0 \pmod{5}.
\end{align}
\end{thm}

The paper is organized as follows. In Section 2, we recall some definitions and properties on modular forms. In Section 3, we give the proof of theorem \ref{b11-p-thm}. In Section 4, we give the proofs of the remain theorems. We conclude the paper in Section 5 by giving some explicit examples.

\section{Preliminaries}
To make this paper self-contained, we recall some definitions and basic knowledge on modular forms. For more details, see, for example \cite{Koblitz-1993, Ono-2004}.

Let $k$ be an integer, $N$ be a positive integer  and $\chi$ be a Nebentypus character. We use $M_k(\Gamma_0(N), \chi)$ to denote
the space of holomorphic modular forms on $\Gamma_0(N)$ of weight $k$ with character $\chi$. The corresponding space of cusp forms
is denoted by $S_k(\Gamma_0(N), \chi)$.
If $\chi$ is the trivial character, we shall write $M_k(\Gamma_0(N))$ and $S_k(\Gamma_0(N))$ for short.
Moreover, we write $SL_2(\mathbb{Z})$ for $\Gamma_0(1)$.

\begin{defi}
Let $m$ be a positive integer and $f(z)=\sum_{n=0}^{\infty}a(n)q^n \in M_k(\Gamma_0(N),\chi)$, then the action of Hecke operator $T_m$ on $f(z)$ is defined by
\begin{equation*}
f(z)\mid T_m
:=\sum_{n=0}^{\infty}\Big(\sum_{d|{\rm gcd}(m,n)}\chi(d)d^{k-1}a(mn/d^2) \Big)q^n.
\end{equation*}
In particular, if $m=p$ is prime, we have
\begin{equation}\label{Tp}
  f(z)\mid T_p:=\sum_{n=0}^{\infty}\Big(a(pn)+\chi(p)p^{k-1}a(n/p)\Big)q^n.
\end{equation}
Note that $a(n)=0$ unless $n$ is a nonnegative integer.
\end{defi}

\begin{defi}
Let $f(z)=\sum_{n=0}^{\infty}a(n)q^n \in M_k(\Gamma_0(N),\chi)$ be a modular form, then $f(z)$ is called a Hecke eigenform if for every $m\ge 2$ there exists a complex number $\lambda(m)$ for which
\begin{equation}\label{H-eig}
  f(z)\mid T_m = \lambda(m)f(z).
\end{equation}
\end{defi}

We recall that Dedekind's eta-function $\eta(z)$ is defined by
\[
\eta(z):=q^{\frac{1}{24}}(q;q)_{\infty},
\]
which is a holomorphic modular form of weight $\frac{1}{2}$ and is non-vanishing on the upper half complex plane $\mathcal{H}$.

\begin{pro}[Gordon, Hughes, Newman]\label{GHN}
If $f(z)=\prod_{\delta |N}\eta^{r_{\delta}}(\delta z)$ is an eta-quotient with $r_{\delta}\in \mathbb{Z}$, $k=\frac{1}{2}\sum_{\delta |N} r_{\delta}\in \mathbb{Z}$, and with the additional properties that
\begin{equation}\label{cond-1}
  \sum_{\delta |N}\delta r_{\delta}\equiv 0 \pmod{24}
\end{equation}
and
\begin{equation}\label{cond-2}
 \sum_{\delta |N}\frac{Nr_{\delta}}{\delta}\equiv 0 \pmod{24},
\end{equation}
then $f(z)$ satisfies
\[
f\Big(\frac{az+b}{cz+d}\Big)=\chi(d)(cz+d)^kf(z)
\]
for all $\Big(
           \begin{array}{cc}
             a & b \\
             c & d \\
           \end{array}
         \Big)
\in \Gamma_0(N)$. Here the character $\chi$ is defined by $\chi(d)=\left(\frac{(-1)^ks}{d} \right)$ and $s=\prod_{\delta|N}\delta^{r^{\delta}}$.
\end{pro}

For verifying that an eta-quotient $f(z)$ is a modular form, we also need the following proposition to compute the orders at cusps.

\begin{pro}\label{cusp}
Let $c,d$ and $N$ be positive integers with $d|N$ and $gcd(c,d)=1$. If $f(z)$ is an eta-quotient satisfying the conditions \eqref{cond-1} and \eqref{cond-2} of Proposition \ref{GHN} for $N$, then the order of vanishing of $f(z)$ at the cusp $\frac{c}{d}$ is
\[
\frac{N}{24d (d,\frac{N}{d})}\sum_{d|N}\frac{(d,\delta)^2r_{\delta}}{\delta}.
\]

\end{pro}

\section{Proof of Theorem \ref{b11-p-thm}}
In this section, we give the proof of Theorem \ref{b11-p-thm}  by using the following \emph{Vanishing Property} given by Hou et al. \cite{Hou-Sun-Zhang-2015}.

\begin{lem}[Vanishing Property] \label{lem-h}
Let $p$ be a prime and
\[
  F(q) = \sum_{n=0}^\infty a(n)q^n = \sum_{k,l=-\infty}^\infty c(k,l)q^{\theta(k,l)}.
\]
Suppose that there exist integers $\theta_0, r,s$ and an invertible transformation $\sigma \colon \mathbb{Z}^2 \to \mathbb{Z}^2$ satisfying the following three conditions
\begin{enumerate}
\item[{\rm (a)}]
the congruence $\theta(k, l) \equiv \theta_0 \pmod p$ has a unique solution $k \equiv r \pmod{p}$ and $l \equiv s \pmod{p}$ in $\mathbb{Z}_p^2$;

\item[{\rm (b)}]
$\theta(pk+r, pl+s) = p^2 \theta(\sigma(k,l)) + \theta_0$;

\item[{\rm (c)}]
$c(pk+r,pl+s) = \lambda(p) \cdot c(\sigma(k,l))$, where $\lambda(p)$ is a constant independent of $k$ and $l$.
\end{enumerate}
Then the following two assertions hold.
\begin{enumerate}
\item[{\rm (1)}] For any integer $n$, we have
\[
  a(p^2 n + \theta_0) = \lambda(p) \cdot a(n).
\]

\item[{\rm (2)}]
For any integer $n$ with $p \nmid n$, we have
\begin{equation}\label{eq:pj}
a(p n + \theta_0) = 0.
\end{equation}

\end{enumerate}
\end{lem}

{\noindent \emph{Proof of Theorem \ref{b11-p-thm}}}. We have
\[
\sum_{n=0}^{\infty}b_{11}(n)q^n
= \frac{(q^{11},q^{11})_{\infty}}{(q;q)_\infty}
\equiv (q;q)_{\infty}^{10}
 \pmod{11}.
\]
Denote by $h(q):=(q;q)_{\infty}^{10}=\sum_{n=0}^\infty a(n)q^n$.
Recall that the Winquist's identity is given by (see \cite{Kang-1997})
\begin{align}\label{win}
 &\sum_{k=-\infty}^{\infty}\sum_{l=-\infty}^{\infty}
 (-1)^{k+l}q^{(3k^2+3l^2+3k+l)/2}\nonumber\\[5pt]
 &\hspace{1.5cm}\times(a^{-3k}b^{-3l}-a^{-3k}b^{3l+1}
      -a^{-3l+1}b^{-3k-1}+a^{3l+2}b^{-3k-1})\nonumber\\[5pt]
 &\hspace{2.5cm}= (a,a^{-1}q,b,b^{-1}q,ab,a^{-1}b^{-1}q,ab^{-1},a^{-1}bq;q)_{\infty}(q;q)_{\infty}^2,
\end{align}
where $q$ is any complex number with $|q| < 1$ and $a, b$ are any nonzero complex numbers.
By elementary manipulations, we can see that identity \eqref{win} can be rewritten as
\begin{align}\label{win-1}
 &\sum_{k=0}^{\infty}\sum_{l=-\infty}^{\infty}
 (-1)^{k+l}q^{(3k^2+3l^2+3k+l)/2}\nonumber\\[5pt]
 &\hspace{1.5cm}\times(a^{-3k}b^{-3l}-a^{-3k}b^{3l+1}
     -a^{3k+3}b^{-3l}+a^{3k+3}b^{3l+1}\nonumber\\[5pt]
 &\hspace{2.5cm}  +a^{-3l+1}b^{3k+2}-a^{3l+2}b^{3k+2}
      -a^{-3l+1}b^{-3k-1}+a^{3l+2}b^{-3k-1})\nonumber\\[5pt]
 &\hspace{3.5cm}= (a,a^{-1}q,b,b^{-1}q,ab,a^{-1}b^{-1}q,ab^{-1},a^{-1}bq;q)_{\infty}
 (q;q)_{\infty}^2.
\end{align}
Taking the limits $a\rightarrow 1$ and $b\rightarrow 1$ in \eqref{win-1} (by using L'Hospital Rule), we get
\begin{align}\label{b11-c}
 (q;q)_{\infty}^{10}
& =\sum_{k=0}^{\infty}\sum_{l=-\infty}^{\infty}
 (-1)^{k+l}(2k+1)(6l+1)\nonumber\\[5pt]
&\hspace{2cm} \times\Big(\frac{(3k+1)(3k+2)}{2}-\frac{(3l)(3l+1)}{2}\Big)
 q^{\frac{3k(k+1)}{2}+\frac{l(3l+1)}{2}}.
\end{align}
Assume that
\[
(q;q)_{\infty}^{10}
=\sum_{k=0}^{\infty}\sum_{l=-\infty}^{\infty}c(k,l)q^{\theta(k,l)},
\]
where
\[
c(k,l)=(-1)^{k+l}(2k+1)(6l+1)\Big(\frac{(3k+1)(3k+2)}{2}-\frac{(3l)(3l+1)}{2}\Big),
\]
and
\[
\theta(k,l)=\frac{3k(k+1)}{2}+\frac{l(3l+1)}{2}.
\]
Notice that
\[
\theta(k,l) = \frac{3}{2} \left( \left( k+\frac{1}{2} \right)^2 + \left( l+\frac{1}{6} \right)^2 \right) - \frac{5}{12}.
\]
For any $1 \le i \le \alpha+1$ and prime $p_i\ge 5$, we have
\[
  \theta(k,l) \equiv - \frac{5}{12} \pmod{p_i}   \quad \Leftrightarrow \quad  \left( k+\frac{1}{2} \right)^2 + \left( l+\frac{1}{6} \right)^2 \equiv 0 \pmod{p_i}
\]
Since $p_i \equiv 3 \pmod 4$, $-1$ is not a quadratic residue modulo $p_i$. Hence
\begin{align*}
\left( k+\frac{1}{2} \right)^2 \equiv - \left( l+\frac{1}{6} \right)^2 \pmod{p_i} \quad \Leftrightarrow \quad k \equiv \frac{p_i-1}{2} \quad \& \quad l \equiv \frac{\pm p_i-1}{6} \pmod{p_i},
\end{align*}
where $\pm$ makes sure that $(\pm p_i-1)/6$ should be an integer.
Furthermore, we have
\[
  \theta \left( kp_i+\frac{p_i-1}{2}, lp_i+\frac{\pm p_i-1}{6} \right) = p_i^2 \theta(k,\pm l) + \frac{5(p_i^2-1)}{12}
\]
and
\[
  c \left( kp_i+\frac{p_i-1}{2}, lp_i+\frac{\pm p_i-1}{6} \right) = p_i^4 c(k,\pm l).
\]
We thus obtain the following recurrence relation from Lemma~\ref{lem-h} (1).
\begin{equation}\label{eq:rec}
  a \left( p_i^2 n + \frac{5(p_i^2-1)}{12} \right) = p_i^4a(n).
\end{equation}
Iteratively using recursion \eqref{eq:rec}, we obtain that
\begin{align*}
a\left(p_1^2\cdots p_{\alpha}^2p_{\alpha +1}n
   +\frac{5(p_1^2\cdots p_{\alpha+1}^2-1)}{12}\right)
 = p_1^4\cdots p_{\alpha}^4 \cdot a\left(p_{\alpha +1}n
   +\frac{5(p_{\alpha+1}^2-1)}{12}\right) .
\end{align*}
By Lemma~\ref{lem-h} (2), $a\left(p_{\alpha +1}n
   +\frac{5(p_{\alpha+1}^2-1)}{12}\right) \not=0$ only when $p_{\alpha+1} \mid n$. Therefore,
\[
  \sum_{n=0}^\infty a\left(p_{\alpha +1}n
   +\frac{5(p_{\alpha+1}^2-1)}{12}\right) q^n = \sum_{n'=0}^\infty a\left(p_{\alpha +1}^2 n'
   +\frac{5(p_{\alpha+1}^2-1)}{12}\right) q^{p_{\alpha+1}n'}.
\]
Using recursion \eqref{eq:rec} once again, the above sum reduces to
\[
  \sum_{n'=0}^\infty a(n') q^{p_{\alpha+1}n'} = h(q^{p_{\alpha+1}})^2.
\]
Thus we obtain that
\begin{equation*}
  \sum_{n=0}^\infty a \left(p_1^2\cdots p_{\alpha}^2p_{\alpha +1}n
   +\frac{5(p_1^2\cdots p_{\alpha+1}^2-1)}{12}\right)q^n
  =  p_1^4\cdots p_{\alpha+1}^4 \cdot h(q^{p_{\alpha +1}})^2.
\end{equation*}
This gives
\begin{equation}\label{11-r}
  \sum_{n=0}^\infty b_{11} \left(p_1^2\cdots p_{\alpha}^2p_{\alpha +1}n
   +\frac{5(p_1^2\cdots p_{\alpha+1}^2-1)}{12}\right)q^n
  \equiv  p_1^4\cdots p_{\alpha+1}^4 \cdot h(q^{p_{\alpha +1}})^2 \pmod{11},
\end{equation}
which finishes the proof of Theorem \ref{b11-p-thm} by replacing $n$ by $p_{\alpha +1}n+j$ in \eqref{11-r}.
\qed

\section{Proofs of Other Theorems}
In this section, we prove the remaining theorems by using the theory of Hecke eigenforms.

{\noindent \emph{Proof of Theorem \ref{b3-p-thm}}}. We have
\[
\sum_{n=0}^{\infty}b_3(n)q^n
=\frac{(q^3;q^3)_{\infty}}{(q;q)_{\infty}}
\equiv (q;q)_{\infty}^2 \pmod{3}.
\]
Replacing $q$ by $q^{12}$ and multiplying $q$ on both sides of the above identity, we deduce that
\[
\sum_{n=0}^{\infty}b_3(n)q^{12n+1}
\equiv q(q^{12};q^{12})_{\infty}^2 = \eta^2(12z) \pmod{3}.
\]
Let $\eta^2(12z)=\sum\limits_{n=1}^{\infty}a(n)q^n$. It is clear that $a(n)=0$ if $n\not\equiv 1 \pmod{12}$. And also, for all $n\geq 0$,
\begin{equation}\label{a-b3}
b_3(n)\equiv a(12n+1) \pmod{3}.
\end{equation}
Notice that $\eta^2(12z) \in S_1\left(\Gamma_0(144),(\frac{-1}{d})\right)$ is an eigenform (see, for example \cite{Vos-2012}).
By \eqref{Tp} and \eqref{H-eig}, we have
\begin{align*}
  \eta^2(12z)\mid T_p
= \sum_{n=1}^{\infty}\left(a(pn)+\Big(\frac{-1}{p}\Big)a(n/p) \right)q^n
= \lambda(p) \sum_{n=1}^{\infty}a(n)q^n,
\end{align*}
which implies
\begin{align}\label{b3-eig}
a(pn)+\Big(\frac{-1}{p}\Big)a(n/p)=\lambda(p)a(n).
\end{align}
By setting $n=1$ in \eqref{b3-eig} and noting that $a(1)=1$, we can see $a(p)=\lambda(p)$.
Since $p\not\equiv 1 \pmod{12}$, we thus obtain
\begin{align}\label{b3-ap-0}
\lambda(p) = a(p)=0.
\end{align}
Combining \eqref{b3-eig} with \eqref{b3-ap-0}, we derive that
\begin{align}\label{b3-a-re}
 a(pn)=-\Big(\frac{-1}{p}\Big)a(n/p).
\end{align}
From \eqref{b3-a-re}, we obtain that for all $n\ge 0$ and $p\nmid r$,
\begin{align}\label{b3-a-0}
 a(p^2n+pr)=0,
\end{align}
and
\begin{align}\label{b3-a-1}
 a(p^2n)=-\Big(\frac{-1}{p}\Big)a(n).
\end{align}
Substituting $n$ by $12n-pr+1$ in \eqref{b3-a-0} and together with \eqref{a-b3}, we find that
\begin{equation}\label{b3-0}
 b_3\Big(p^2n+\frac{p^2-1}{12}+p\frac{1-p^2}{12}r\Big)
\equiv 0 \pmod{3}.
\end{equation}
Since $p\ge 5$ is prime, we have $12 \mid (1-p^2) $, and ${\rm gcd}(\frac{1-p^2}{12},p)=1$. Hence when $r$ runs over a residue system excluding the multiple of $p$, so does $\frac{1-p^2}{12}r$. Thus \eqref{b3-0} can be rewritten as
\begin{equation}\label{b3-id}
  b_3\Big(p^2n+\frac{p^2-1}{12}+pj \Big)
\equiv 0 \pmod{3},
\end{equation}
where $p\nmid j$.
Replacing $n$ by $12n+1$ in \eqref{b3-a-1} and together with \eqref{a-b3}, we have
\begin{equation}\label{b3-re}
 b_3\Big(p^2n+\frac{p^2-1}{12}\Big)
\equiv -\Big(\frac{-1}{p}\Big) b_3(n) \pmod{3}.
\end{equation}
If  $p_i\geq 5$ are primes such that $p_i\not\equiv 1 \pmod{12}$, then
iteratively using recursion \eqref{b3-re}, we obtain that
\begin{align}\label{b3-p-re}
 b_3\left(p_1^2\cdots p_{\alpha}^2n+\frac{p_1^2\cdots p_{\alpha}^2-1}{12}\right)
 \equiv (-1)^{\alpha}\Big(\frac{-1}{p_1}\Big)\cdots\Big(\frac{-1}{p_{\alpha}}\Big)b_3(n) \pmod{3}.
\end{align}
Replacing $p$ by $p_{\alpha+1}$ in \eqref{b3-id} and combining with \eqref{b3-p-re}, we derive that
\begin{align*}
 b_3\left(p_1^2\cdots p_{\alpha+1}^2n+\frac{p_1^2\cdots p_{\alpha}^2p_{\alpha +1}(12j+p_{\alpha +1})-1}{12}\right)\equiv 0 \pmod{3}.
\end{align*}
This finishes the proof of Theorem \ref{b3-p-thm}.
\qed

Remark that when $p_i\ge 5$ are prime and $p_i\equiv 3 \pmod{4}$, the congruences were proven by Hou, Sun and Zhang \cite{Hou-Sun-Zhang-2015} by using binary quadratic forms. It is worthy noticing that we obtain a much more family of congruences  for more primes $p_i\not\equiv 1 \pmod{12}$ by using Hecke operator on eigenform.

{\noindent \emph{Proof of Theorem \ref{b13-p-thm}}}. We have
\[
\sum_{n=0}^{\infty}b_{13}(n)q^n
=\frac{(q^{13};q^{13})_{\infty}}{(q;q)_{\infty}}
\equiv (q;q)_{\infty}^{12} \pmod{13}.
\]
Replacing $q$ by $q^2$ and then multiplying $q$ both sides of the above identity, we deduce that
\[
\sum_{n=0}^{\infty}b_{13}(n)q^{2n+1}
\equiv q(q^2;q^2)_{\infty}^{12}
=\eta^{12}(2z) \pmod{13}.
\]
Denote by $\eta^{12}(2z)=\sum_{n=1}^{\infty}a(n)q^n$, we then have $a(n)=0$ if $n$ is even and
\begin{equation}\label{b13-a}
  b_{13}(n)\equiv a(2n+1) \pmod{13}.
\end{equation}
Notice that $\eta^{12}(2z) \in S_6(\Gamma_0(4))$ is an eigenform (see, for example, \cite{Vos-2012}). Then by \eqref{Tp} and \eqref{H-eig}, we have that
\[
\eta^{12}(2z)\mid T_p = \sum_{n=1}^{\infty}\Big(a(pn)+p^5a(n/p)\Big)q^n
=\lambda(p) \sum_{n=1}^{\infty}a(n)q^n,
\]
which implies
\begin{equation}\label{13-a-eig}
 a(pn)+p^5a(n/p)=\lambda(p)a(n).
\end{equation}
Setting $n=1$ in \eqref{13-a-eig} and since $a(1)=1$, we find
$a(p)=\lambda(p)$.
Thus \eqref{13-a-eig} becomes
\begin{equation}\label{13-a-1}
  a(p^2n)+p^5a(n)=a(p)a(pn),
\end{equation}
and for $p\nmid r$,
\begin{equation}\label{13-a-2}
  a(p^2n+pr)=a(p)a(pn+r).
\end{equation}
By \eqref{b13-a} and \eqref{13-a-1}, we derive that
\begin{equation}\label{b13-c-1}
 b_{13}\Big(p^2n+\frac{p^2-1}{2}\Big)+p^5b_{13}(n)
 \equiv b_{13}\Big(\frac{p-1}{2}\Big)b_{13}\Big(pn+\frac{p-1}{2}\Big) \pmod{13}.
\end{equation}
If $ b_{13}(\frac{p-1}{2})\equiv 0 \pmod{13}$, \eqref{b13-c-1} implies that
\begin{equation}\label{b13-r-1}
 b_{13}\Big(p^2n+\frac{p^2-1}{2}\Big) \equiv -p^5b_{13}(n) \pmod{13}.
\end{equation}
Replacing $n$ by $2n-pr+1$ in \eqref{13-a-2} and together with \eqref{b13-a}, we obtain that
\begin{equation}\label{b13-c-2}
 b_{13}\Big(p^2n+\frac{p^2-1}{2}+p\frac{1-p^2}{2}r\Big)
 \equiv b_{13}\Big(\frac{p-1}{2}\Big)b_{13}\Big(pn+\frac{p-1}{2}+\frac{1-p^2}{2}r\Big) \pmod{13}.
\end{equation}
Since $p$ is odd prime, so $2 \mid (1-p^2) $, and ${\rm gcd}(\frac{1-p^2}{2},p)=1$. Hence when $r$ runs over a residue system excluding the multiple of $p$, so does $\frac{1-p^2}{2}r$. Thus \eqref{b13-c-2} can be rewritten as
\begin{equation}\label{b13-r-2}
 b_{13}\Big(p^2n+\frac{p^2-1}{2}+pj\Big)
 \equiv b_{13}\Big(\frac{p-1}{2}\Big)b_{13}\Big(pn+\frac{p-1}{2}+j\Big) \pmod{13},
\end{equation}
where $p\nmid j$.
If $ b_{13}(\frac{p-1}{2})\equiv 0 \pmod{13}$, \eqref{b13-r-2} implies that
\begin{equation}\label{b13-r-0}
 b_{13}\Big(p^2n+\frac{p^2-1}{2}+pj\Big) \equiv 0 \pmod{13},
\end{equation}
where $p\nmid j$.
Therefore, for $1\leq i \le \alpha+1$, if $ b_{13}(\frac{p_i-1}{2})\equiv 0 \pmod{13}$, then by using the recursion \eqref{b13-r-1} iteratively , we obtain
\begin{equation}\label{b13-pi}
 b_{13}\Big(p_1^2\cdots p_{\alpha}^2n+\frac{p_1^2\cdots p_{\alpha}^2-1}{2}\Big) \equiv (-1)^{\alpha}p_1^5\cdots p_{\alpha}^5 b_{13}(n) \pmod{13}.
\end{equation}
Furthermore, by \eqref{b13-r-0}, we get
\begin{equation}\label{b13-p-0}
 b_{13}\Big( p_{\alpha+1}^2n+\frac{p_{\alpha+1}^2-1}{2}+p_{\alpha+1}j\Big) \equiv 0 \pmod{13}.
\end{equation}
Combining \eqref{b13-pi} with \eqref{b13-p-0}, we complete the proof of Theorem \ref{b13-p-thm}.
\qed

{\noindent \emph{Proof of Theorem \ref{b25-p-thm}}}. We have
\[
\sum_{n=0}^{\infty}b_{25}(n)q^n
=\frac{(q^{25};q^{25})_{\infty}}{(q;q)_{\infty}}
\equiv (q;q)_{\infty}^{24} \pmod{5}.
\]
Multiplying $q$ on both sides of the above identity, we deduce that
\[
\sum_{n=0}^{\infty}b_{25}(n)q^{n+1}
\equiv q(q;q)_{\infty}^{24}
=\eta^{24}(z) \pmod{5}.
\]
Denote by $\eta^{24}(z)=\sum_{n=1}^{\infty}a(n)q^n$, we then have
\begin{equation}\label{b25-a}
  b_{25}(n)\equiv a(n+1) \pmod{5}.
\end{equation}
Notice that $\eta^{24}(z) \in S_{12}(SL_2(\mathbb{Z}))$ is an eigenform (see, for example, \cite{Vos-2012}). Then by \eqref{Tp} and \eqref{H-eig}, we have that
\[
\eta^{24}(z)\mid T_p = \sum_{n=1}^{\infty}\Big(a(pn)+p^{11}a(n/p)\Big)q^n
=\lambda(p) \sum_{n=1}^{\infty}a(n)q^n,
\]
which implies
\begin{equation}\label{25-a-eig}
 a(pn)+p^{11}a(n/p)=\lambda(p)a(n).
\end{equation}
Setting $n=1$ in \eqref{25-a-eig} and since $a(1)=1$, we follow that $a(p)=\lambda(p)$. Thus \eqref{25-a-eig} becomes
\begin{equation}\label{25-a-1}
  a(p^2n)+p^{11}a(n)=a(p)a(pn),
\end{equation}
and for $p\nmid j$,
\begin{equation}\label{25-a-2}
  a(p^2n+pr)=a(p)a(pn+j).
\end{equation}
If $b_{25}(p-1)\equiv 0 \pmod{5}$, we obtain from \eqref{b25-a} and \eqref{25-a-1}  that
\begin{equation}\label{25-b-re-1}
b_{25}(p^2n+p^2-1)\equiv -p^{11}b_{25}(n) \pmod{5}.
\end{equation}
Also, by \eqref{b25-a} and \eqref{25-a-2}, we obtain that for $p\nmid j$,
\begin{equation}\label{25-b-re-2}
b_{25}(p^2n+pj+p^2-1)\equiv 0 \pmod{5}.
\end{equation}
Since $p_i$ are primes and $b_{25}(p_i-1)\equiv 0 \pmod{5}$, by using recursion \eqref{25-b-re-1} iteratively, we get
\begin{equation}\label{25-c-1}
  b_{25}(p_1^2\cdots p_{\alpha}^2n+p_1^2\cdots p_{\alpha}^2-1)
 \equiv (-1)^{\alpha}p_1^{11}\cdots p_{\alpha}^{11}b_{25}(n) \pmod{5}.
\end{equation}
Further, by \eqref{25-b-re-2}, we obtain
\begin{equation}\label{25-c-2}
b_{25}(p_{\alpha+1}^2n+p_{\alpha+1}j+p_{\alpha+1}^2-1)\equiv 0 \pmod{5}.
\end{equation}
Together with \eqref{25-c-1} and \eqref{25-c-2}, we complete the proof of Theorem \ref{b25-p-thm}.
\qed

\section{Some Examples}
Now we give some explicit congruences from our theorems in the previous sections to conclude this paper.

Let $\alpha$ be a positive integer, $p$ be a prime and $j$ be an integer with $p \nmid j$.
By setting all the primes $p_1, p_2,\ldots, p_{\alpha+1}$ to be equal to the same prime $p$,
we can derive the following infinite families of congruences for $b_\ell(n)$. Note that
Carlson and Webb \cite{Carlson-Webb-2014} have obtained some similar congruences of $b_\ell(n)$ for
$\ell\in\{10,15,20\}$.

\begin{enumerate}
\item
If $p\geq 5$ and $p \not\equiv 1 \pmod{12}$, we have
\begin{align*}
b_3 \left( p^{2\alpha}n+p^{2\alpha-1}j+\frac{p^{2\alpha}-1}{12} \right)\equiv 0 \pmod{3}.
\end{align*}

\item
If $p\geq 5$ and $p \not\equiv 3 \pmod{4}$, we have
\[
b_{11} \left( p^{2\alpha}n+p^{2\alpha-1}j + \frac{5(p^{2\alpha}-1)}{12} \right)\equiv 0 \pmod{11}.
\]

\item
If $p\geq 3$ such that $b_{13}(\frac{p-1}{2})\equiv 0 \pmod{13}$, we have
\[
b_{13} \left( p^{2\alpha}n+p^{2\alpha-1}j + \frac{p^{2\alpha}-1}{2} \right)\equiv 0 \pmod{13}.
\]

\item
If  $b_{25}(p-1)\equiv 0 \pmod{5}$, we have
\[
b_{25} \left( p^{2\alpha}n+p^{2\alpha-1}j + p^{2\alpha}-1 \right)\equiv 0 \pmod{5}.
\]
\end{enumerate}

Now setting $\alpha=1$ and taking some explicit primes in the above congruence relations, we obtain the following congruences.
\begin{enumerate}
\item
For $n\geq 0$ and $j\not\equiv 1 \pmod{5}$, we have
\begin{align*}
b_3(25n+5j-3)\equiv 0 \pmod{3}.
\end{align*}

\item
For $n\geq 0$ and $j\not\equiv 3 \pmod{7}$, we have
\begin{align*}
b_{11}(49n+7j-1)\equiv 0 \pmod{11}.
\end{align*}

\item We find $b_{13}(75)\equiv 0 \pmod{13}$ by Maple. Hence we have
\begin{align*}
b_{13}(22801n+151j+75)\equiv 0 \pmod{13},
\end{align*}
for $n\geq 0$ and $j\not\equiv 75 \pmod{151}$.

\item We find $b_{25}(5-1)\equiv 0 \pmod{5}$ by Maple. Hence we have
\begin{align*}
b_{25}(25n+5j-1)\equiv 0 \pmod{25},
\end{align*}
for $n\geq 0$ and $j\not\equiv 0 \pmod{5}$.
\end{enumerate}

\vspace{.3cm}

{\noindent \bf Acknowledgments.}
We thank Professor Qing-Hu Hou for his helpful suggestions and discussions.
This work was supported by the PCSIRT Project
of the Ministry of Education and the National Science Foundation of China. The first author was also supported by
the National Science Foundation of China (Tianyuan Fund for Mathematics, No.~11426166) and the Project Sponsored
by the Scientific Research Foundation of Tianjin University of Technology and Education.

\end{document}